\documentclass[11pt,a4paper]{article}

\usepackage{latexsym}
\usepackage[dvips]{graphicx}
\graphicspath{{Figures/}}

\usepackage{a4wide,amsthm,amssymb,amsfonts,color,amsmath,bbm}

\newtheorem{theorem}{Theorem}[section]
\newtheorem{lemma}[theorem]{Lemma}
\newtheorem{corollary}[theorem]{Corollary}
\newtheorem{proposition}[theorem]{Proposition}

\newtheorem{definition}[theorem]{Definition}
\newtheorem{observation}[theorem]{Observation}

\def\ITEMMACRO #1 ??? #2 ???{\par\vskip4pt\noindent%
\hangindent=#2em\setbox0\hbox{#1\kern4pt}%
\ifdim\wd0<\hangindent\setbox0\hbox to\hangindent{\hss#1\kern7pt}\fi%
\box0\ignorespaces}

\def\Item(#1){\ITEMMACRO {\rm (#1)} ??? 1.8 ???}

\let\Bitem=\bItem
\def\BrackItem[#1]{\ITEMMACRO [#1] ??? 1.8 ???}

\def\CG{\widetilde{G}}

\def\ni{\noindent}
\def\Proof{\ni{\sl Proof.}\ }

\def\Claim#1.{\medbreak\ni{\bf Claim~#1.}\ }
\def\qedclaim{\hfill$\triangle$\smallskip}

\def\NN{\hbox{\sf I\kern-1ptN}}

\def\interior{\operatorname{int}}

\begin{document}

\title{Binary Labelings for Plane Quadrangulations
and their Relatives
}

\author{
      \normalsize Stefan Felsner\\
      \small\sf Institut f\"ur Mathematik,\\[-1mm]
      \small\sf Technische Universit\"at Berlin.\\[-1mm]
      \small\sf {\tt felsner@math.tu-berlin.de}
\and
      \normalsize Clemens Huemer
      \footnote{Partially supported by projects MEC MTM2006-01267 and Gen. Cat. 2005SGR00692.}\\
      \small\sf Departament de Matem\`{a}tica Aplicada II,\\[-1mm]
      \small\sf Universitat Polit\`{e}cnica de Catalunya.\\[-1mm]
      \small\sf {\tt clemens.huemer@upc.edu}
\and
      \normalsize Sarah Kappes
      \footnote{Supported by the Deutsche
       Forschungsgemeinschaft trough the international research training group
      `Combinatorics, Geometry, and Computation' (No. GRK 588/2).}\\
      \small\sf Institut f\"ur Mathematik,\\[-1mm]
      \small\sf Technische Universit\"at Berlin.\\[-1mm]
      \small\sf {\tt kappes@math.tu-berlin.de}
\and
      \normalsize David Orden
      \footnote{Research partially supported by grants MTM2005-08618-C02-02 and S-0505/DPI/0235-02.}\\
      \small\sf Departamento de Matem\'aticas,\\[-1mm]
      \small\sf Universidad de Alcal\'a. \\[-1mm]
      \small\sf {\tt david.orden@uah.es}
}

\date{\vbox{\vskip-35mm}}
\maketitle

\begin{abstract}
Motivated by the bijection between Schnyder labelings of a plane triangulation
and partitions of its inner edges into three trees, we look for binary
labelings for quadrangulations (whose edges can be partitioned into two trees).
Our labeling resembles many of the properties of
Schnyder's one for triangulations: Apart from being in bijection with
tree decompositions, paths in these trees allow to define the regions of a
vertex such that counting faces in them yields an algorithm for embedding the
quadrangulation, in this case on a~2-book. Furthermore, as
Schnyder labelings have been extended to 3-connected plane graphs,
we are able to extend our labeling from quadrangulations to a larger
class of 2-connected bipartite graphs.
Finally, we propose a binary labeling for Laman graphs.
\end{abstract}



\section{Introduction}\label{sec:introduction}

Schnyder labelings are by now a classical tool to deal with planar
graphs. A \emph{Schnyder labeling} is a special labeling of the angles
of a plane graph with three colors. Schnyder~\cite{Schnyder}
introduced this concept for triangulations, i.e., maximal (in the
number of edges) planar graphs. He showed that these angle labelings
are in bijection with \emph{Schnyder woods}, i.e., special partitions
of the inner edges of the triangulation into three trees.  A main
application of Schnyder woods are straight-line embeddings of
triangulations on small grids.  Felsner~\cite{Felsner}
generalized the concepts of Schnyder labelings and Schnyder
woods to the larger class of $3$-connected plane graphs.

The present work is motivated by the fact that
{\emph{quadrangulations}}, i.e., maximal bipartite planar graphs,
admit a decomposition of the edge set
into two trees. Our aim is to look for a closer resemblance of
Schnyder structures in these cases. In particular we
study angle labelings with two colors.

In Section~\ref{sec:binary_labelings} we define \emph{weak labelings}
for plane graphs. A weak labeling induces a 2-coloring and a
2-orientation of the $2n-4$ edges (being~$n$ the number of vertices). We show that weak labelings are
indeed in bijection with a pair of 2-orientations, one for the graph
and another for an appropriately defined dual. This allows the
characterization and efficient recognition of graphs admitting a weak
labeling.

In Section~\ref{ssec:strong-quadrangulations} we define \emph{strong
labelings} as a subclass of weak labelings.  A graph admitting a
strong labeling has to be a plane quadrangulation.  We show that
strong labelings indeed resemble many properties of Schnyder
labelings:

\begin{itemize}

\item
  Strong labelings induce a partition of the quadrangulation into two
  oriented trees with Schnyder-like properties, see
  Subsection~\ref{ssec:schnyder-like_properties}.
  For the existence of a 2-tree decomposition of a plane quadrangulation
  there are many references, e.g. \cite{AAGHHHKRV,Fraysseix,FMP,FMR95,
   Llado,NW,Petrovic,Ringel}.  The tree
  decomposition induced by the strong labeling has the nice property
  that at each vertex the two trees are ``separated'', i.e,
  around a vertex the edges of each tree appear consecutively.
  Such separating tree decompositions have been previously
  studied in~\cite{Fraysseix}.

\item Strong labelings also allow to obtain an embedding of a
  quadrangulation on a \emph{2-book}, i.e., a mapping of the nodes to
  a line and a non-crossing embedding of the edges in the half-planes
  separated by that line. In our case each halfplane contains the
  edges of one of the two trees. Let $v_1,\dots,v_n$ be the nodes
  ordered along the line.  Then the trees on each page are
  \emph{alternating}, i.e., there are no two edges~$v_iv_j$
  and~$v_jv_k$ with~$i<j<k$. Book embeddings of graphs are well
  studied and have several applications; see e.g.~\cite{Wood}.  For
  the particular case of quadrangulations, the existence of a 2-book
  embedding with a tree on each page was shown in~\cite{FMP}.  Our
  Schnyder-like technique allows to obtain the alternating property.
  Non-crossing alternating trees were studied and counted
  in~\cite{Gelfand}. They have also appeared as one-dimensional analogs
  of pseudo-triangulations~\cite{RSS}.

\end{itemize}

In Section~\ref{sec:genstronglab}, we generalize the notion of
strong labeling such that it is no longer restricted to
quadrangulations. The \emph{generalized strong labelings} are still in
bijection to pairs of trees. This is similar to the generalization of
Schnyder structures in~\cite{Felsner}. The class of bipartite graphs
admitting a generalized strong labeling is characterized in
Subsection~\ref{ssec:admit_genstronglab}.

Finally, in Section~\ref{sec:A binary labeling for plane Laman graphs}
we propose a variant of weak labelings for plane Laman graphs, those with~$n$ vertices and~$2n-3$ edges such
that any induced subgraph on~$k$ vertices has at most~$2k-3$ edges. These graphs arise in the context of rigidity theory~\cite{Laman} and they are strongly connected to pseudo-triangulations~\cite{osss-cpt-2007,Streinu}.

\section{Weak labelings}\label{sec:binary_labelings}

\begin{definition}\rm
Let $G$ be a plane graph. A {\emph{weak labeling}} for $G$ is a
mapping from the angles of $G$ to $\{0,1\}$ which satisfies the
following conditions: 

\Item(G0) {\bf{Special vertices:}} {There are two special
  vertices~$s_0$ and~$s_1$ on the outer face
  of~$G$, such that all angles incident to~$s_i$ are labeled~$i$.}

\Item(G1) {{\bf{Vertex rule:}} For each vertex $v \notin \{s_0,s_1\}$, the
  incident labels form a non-empty interval of $1s$ and a non-empty
  interval of $0s$.}

\Item(G2) {{\bf{Edge rule:}} For each edge, the incident labels
  coincide at one endpoint and differ at the other.}

\Item(G3) {{\bf{Face rule:}} For each face (including the outer face),
  its labels form a non-empty interval of $1s$
  and a non-empty interval of $0s$.}
\end{definition}

\begin{figure}[htb]
    \begin{center}
        \includegraphics{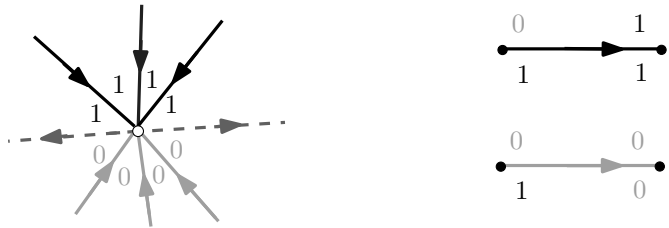}
    \end{center}
    \caption{The orientation induced by a weak labeling, the dashed
      edges may have either color.}
    \label{fig:vertex_rule_edge_rule}
\end{figure}

\begin{observation}\rm
\label{obs:2-coloring_orientation} A weak labeling induces both a
2-coloring and a 2-orientation of the edges: Every edge is colored
according to its endpoint with the two coincident labels and
oriented towards that endpoint. Moreover, the vertex rule implies
that every vertex except~$s_0,s_1$ has outdegree two; such an
orientation will be called \emph{a $2$-orientation}. See
Figures~\ref{fig:vertex_rule_edge_rule} and~\ref{fig:quad_trees}.
\label{obs:weak_2n-4}
It follows that a plane graph with $n$ vertices that admits a weak
labeling must have exactly $2n-4$ edges.
\end{observation}

A quadrangulation on $n$ vertices has $2n-4$ edges and indeed
quadrangulations admit weak labelings (they even admit a stronger
labeling, see Section~\ref{ssec:strong-quadrangulations}). But weak
labelings also exist for some graphs which are not quadrangulations;
consider e.g.  the graph obtained by inserting into the cycle~$C_6$
the edges~$15$ and~$24$. A more complex example is part of
Figure~\ref{fig:X+Xd}.

\subsection{Graphs admitting a weak labeling}
\label{ssec:admit_weak_labeling}

\begin{observation}\label{obs:dual_2-or}\rm
Since the angles of a plane graph $G$ and of the dual $G^*$ are in
bijection we can interpret a weak labeling $\phi$ of $G$ as a labeling
$\phi^*$ of the angles of $G^*$. The edge rule for $\phi$ is the
edge rule for $\phi^*$ and the face rule for $\phi$ implies
the vertex rule for (all!) vertices in the labeling $\phi^*$.
As in the previous observation, we can argue that $\phi^*$ induces an
orientation of the edges of $G^*$ such that every vertex has outdegree
two. One additional property of the orientation induced by $\phi$ can
be noted: The special vertices~$s_0,s_1$ divide the outer face of~$G$
into two arcs~$A_0$ and~$A_1$. Each of these arcs contains a change
of label. Split the dual vertex corresponding to the
outer face of~$G$ into two vertices~$o^*_0$ and $o^*_1$ such that
$o^*_i$ keeps the incidences with the dual edges of $A_i$ and let
$G^*_s$ denote the resulting {\em split-dual} of~$G$.
The orientation induced by~$\phi$ on~$G^*_s$ has outdegree one
at $o^*_0$ and $o^*_1$ and outdegree two at every other vertex;
let us call such an orientation a \emph{a $2^*$-orientation}.
\end{observation}

\begin{proposition}\label{prop:weak_lab+orient}
Let $G$ be a plane graph with special vertices $s_0$ and $s_1$
on the outer face. Weak labelings of $G$ are in bijection to
pairs $(X,X^*)$ where $X$ is a 2-orientation of
$G$ and $X^*$ is a $2^*$-orientation of $G^*_s$.
\end{proposition}

\Proof
The mapping from a weak labeling $\phi$ to a pair $(X_\phi,X^*_\phi)$
of orientations was given in Observations~\ref{obs:2-coloring_orientation}
and~\ref{obs:dual_2-or}.
For the converse construction, we introduce an auxiliary graph:
The {\em completion} $\CG$ of $G$ is obtained by
superimposing $G$ and $G_s^*$ such that exactly the primal-dual pairs of edges cross,
and this crossing is made a new {\em edge-vertex}.
For~$v,e,f$ the numbers of vertices, edges and faces of~$G$,
the completion~$\CG$ of~$G$ has $v+e+f+1$ vertices and $4e$ edges.
The faces of~$\CG$ are (almost) in bijection to the angles of $G$, only
the outer face of~$\CG$  corresponds to two angles,
the outer ones of~$s_0$ and~$s_1$. In order to remedy this,
an exceptional edge~$e_o$ connecting~$o^*_0$ and~$o^*_1$ can be added.

\begin{figure}[htb]
    \begin{center}
        \includegraphics[scale=0.7]{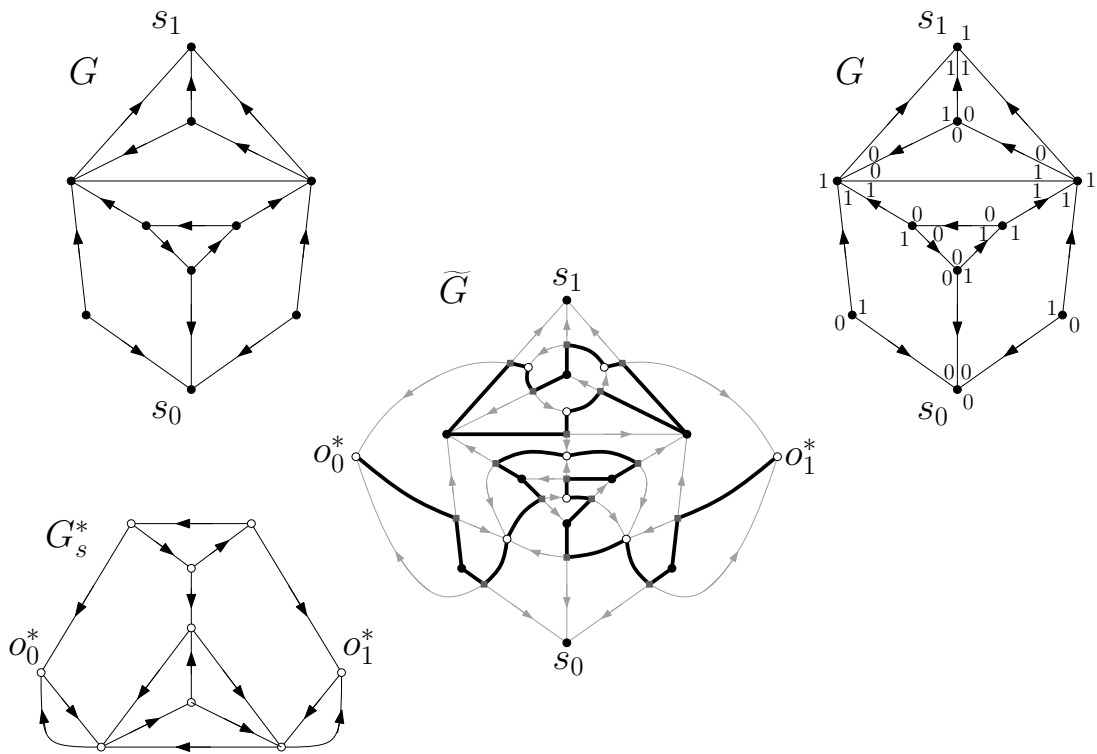}
    \end{center}
    \caption{Orientations $X$ and $X^*$, the relevant edges of~$\CG$ and the resulting
             weak labeling of~$G$.}
    \label{fig:X+Xd}
\end{figure}

Given a 2-orientation $X$ of~$G$ and a $2^*$-orientation $X^*$
of~$G_s^*$, we induce an orientation on~$\CG$ by taking the orientation of
an edge~$e$ for both of its halfedges, see Figure~\ref{fig:X+Xd}. We
seek for a 0-1 coloring of the inner faces of~$\CG$ such that, if~$v_e$
is an edge-vertex and~$a$ an outgoing edge at~$v_e$, then the color of
the two faces incident to~$a$ are the same. we model this by calling~$a$
an {\em irrelevant} edge. Edges that are not irrelevant are {\em relevant}.
Observe that, from the properties of~$X$ and~$X^*$ and
the construction, we obtain:
\Bitem There are no relevant edges
incident to $s_0$ and $s_1$.
\Bitem There is exactly one relevant
edge incident to $o^*_0$ and $o^*_1$.
\Bitem Apart from these
exceptions, every vertex of $\CG$ is incident to exactly two relevant
edges.\smallskip

\ni
It follows that the relevant edges form a union of
disjoint simple cycles and a path from~$o^*_0$ to~$o^*_1$ which, by adding~$e_o$
as relevant edge is also closed into a cycle.
Starting with color~$0$ in the face containing~$s_0$, there is
a unique extension to a 0-1 coloring of the faces in the graph of relevant edges
and hence to a coloring of the faces of~$\CG$. It is routine to check that this
indeed yields a weak labeling and that the two mappings are inverse to each other.
\qed
\medskip

A nice consequence of the above proposition is that it yields a
characterization of plane graphs admitting a weak labeling: Given a
graph $G=(V,E)$ and a function $\alpha: V \to \NN$, an {\em $\alpha$-orientation}
is an orientation $X$ of $G$ such that the outdegree of each $v$ is as
prescribed by $\alpha$, i.e., it is $\alpha(v)$. It is known that $G$
admits an $\alpha$-orientation if $\sum_{v\in V}\alpha(v) = |E|$
and for all $W\subset V$ the number of
edges incident to vertices in $W$ is at least as large as $\alpha(W) =
\sum_{v\in W}\alpha(v)$\label{demand_crit}. Moreover, the question whether $G$ admits an
$\alpha$-orientation can be translated into a flow-problem, hence,
it can be answered in polynomial time. These facts about $\alpha$-orientations
are detailed e.g.~in~\cite{Felsner_latticeStruc}.
A polynomial time recognition algorithm for plane graphs admiting a weak labeling
would first decide whether $G$ and $G_s^*$ admit a $2$- resp.~a $2^*$-orientation.
In the positive case such orientations $X$ and $X^*$ could be transformed into
a weak labeling of $G$. We summarize:

\begin{theorem}
Plane graphs admitting a weak labeling can be recognized in polynomial time.
\end{theorem}

From the theory of $\alpha$-orientations developed
in~\cite{Felsner_latticeStruc} it also follows that the sets of
2-orientations of $G$ and of $2^*$-orientations of $G_s^*$ carry a
natural distributive lattice structure. The product of these two
distributive lattices is a distributive lattice on the set of all weak
labelings of a plane graph. The ideas on how to define the lattice
structure on 2-orientations of a plane graph are explained in
Subsection~\ref{subsec:flips}.

\subsection{Schnyder-like properties for weak labelings}
\label{ssec:schnyder-like_properties_weak}

The orientation and coloring of the edges of a graph induced by a weak
labeling have another interesting property: Let~$G$ be a plane graph
with a weak labeling and let~$T_0$ and~$T_1$ be the edges of colors~0
and~1. Since the edges are oriented according to a 2-orientation, we can
define~$T_i^{-1}$ as the set of edges colored $i$ with their
orientation reversed. In figures,
e.g.~Figure~\ref{fig:vertex_rule_edge_rule}, and sometimes in the text
we will identify color~$0$ with gray and color~$1$ with black.

The following proposition is very much like Schnyder's main lemma
in~\cite{Schnyder89}. However, in a weak labeling a vertex can have
out-degree two in $T_i$ wherefore $T_i$ need not be a tree,
see e.g.~Figure~\ref{fig:X+Xd}.

\begin{proposition}
\label{lem:no_directed_cycles-weak}
If $G$ is a plane graph with a weak labeling, then
there is no directed cycle in~$T_0\cup T_1^{-1}$, nor in $T_1 \cup T_0^{-1}$.
\end{proposition}

\Proof Suppose that there is a directed cycle $C$ in $T_0\cup
T_1^{-1}$. Clearly we may assume that $C$ is simple, hence, has a
well defined interior. Consider $G$ with the original 2-orientation
$T_0\cup T_1$ and define the following counters:
\begin{eqnarray*}
k &=& \#\text{vertices on $C$.}\\
t &=& \#\text{vertices in the interior of $C$.}\\
s &=& \#\text{faces in the interior of $C$.}\\
p &=& \#\text{edges pointing in $T_0\cup T_1$ from $C$ into the interior.}\\
q &=& \#\text{edges on $C$ with two different labels on the inner side.}\\
\end{eqnarray*}
\begin{figure}[htb]
    \begin{center}
        \includegraphics{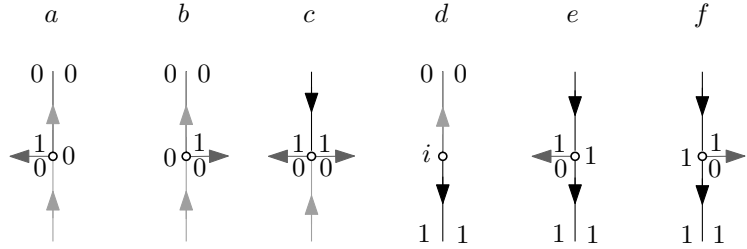}
    \end{center}
    \caption{Vertex types on $C$, where the vertical edges are those on~$C$ and the interior of~$C$ is assumed to be
    to the left of them. The schematic drawings show only the relevant edges.}
    \label{fig:v-types}
\end{figure}

\hbox{\Claim A. $p=q$.}

\ni Figure~\ref{fig:v-types} shows all types of vertices which can
occur on $C$. Associate each edge that has two different labels on
the inner side with its tail-vertex. We find that $q$ equals the
number of
vertices of types $a$, $d$ and~$e$.  The value of $p$ is the number of
vertices of types $a$, $c$ and~$e$.  Since vertices of type $d$
correspond to a transition from $0$-colored edges to $1$-colored edges
while vertices of type $c$ correspond to a transition from $1$-colored
to $0$-colored edges, it follows that they are equinumerous. This
proves the claim. \qedclaim

We now observe that the number~$g$ of edges which are on $C$ or in
the interior of $C$ can be expressed in several ways:
\begin{eqnarray}
g &=& (k+t) + (s+1) -2\\
g &=& 2t + k + p\\
g &=& 2s + k -q
\end{eqnarray}

Formula (1) is nothing but Euler's formula for the graph restricted
to $C$ and its interior. Formula (2) is obtained by counting the
out-degrees: Every vertex in the interior of $C$ has out-degree 2
and the sum of all out-degrees of vertices on $C$ is $k+p$. Formula
(3) follows from counting changes of labels along edges: By the edge
rule (G2) the number of these changes equals the number of edges. By
the face rule (G3) each of the $s$ faces interior to $C$ contributes
two such changes. In addition there are $k-q$ edges on $C$ which
have the label change in the outside.

Subtracting (2) and (3) from the double of (1) yields $0 = -2 -p
+q$, which is a contradiction to Claim~A. \qed

\section{Strong labelings}
\label{ssec:strong-quadrangulations}

Consider weak labelings of a plane graph obeying the following
strong face rule:

\Item(G3$^+$){\bf{Strong face rule:}}
    Each face has exactly one pair of adjacent $0$-labels and one
    pair of adjacent $1$-labels.
    In addition, the edge on the outer face $F_{\sf out}$ which contains $s_0$
    and which has $F_{\sf out}$ to its right when traversed from
    its white end to the black end has two adjacent labels $0$ in $F_{\sf out}$.
\smallskip

\ni
In Observation~\ref{obs:weak_2n-4} we noticed that graphs~$G$ with~$n$ vertices
having a weak labeling must have~$2n-4$ edges. Then, simple counting shows that further
requiring the strong face rule implies:
\begin{observation}\rm
\label{obs:strong_quadrangulation}
Every plane graph that admits a strong
labeling is a quadrangulation, i.e., a maximal bipartite plane graph.
\end{observation}

This allows to state the last of the defining properties of a strong labeling
in a more convenient way:

\begin{definition}\rm
Let $G$ be a quadrangulation with color classes of black and white
vertices. A {\em strong labeling} of $G$ is a
mapping of the angles of $G$ to $\{0,1\}$ which satisfies:

\Item(G0) {\bf{Special vertices:}} {The two black vertices
  on the outer face are named~$s_0$ and~$s_1$.
  All angles incident to~$s_i$ are labeled~$i$.}

\Item(G1) {{\bf{Vertex rule:}} For each vertex $v \notin \{s_0,s_1\}$, the
  incident labels form a non-empty interval of $1s$ and a non-empty
  interval of $0s$.}

\Item(G2) {{\bf{Edge rule:}} For each edge, the incident labels
  coincide at one endpoint and differ at the other.}

\Item(G3$^+_Q$) {{\bf{Strong face rule for quadrangulations:}}
    The labels in each face are 0011 cyclically.
    Reading the labels of the outer face in clockwise order
    starting at $s_0$, they are also 0011.}
\end{definition}

Not every weak labeling of a quadrangulation is strong. For example,
exchanging all labels except those at $s_0$ and $s_1$ turns a strong
labeling into a weak one.

\begin{lemma}[Walking rule]\label{lem:walking}
In a strong labeling of a quadrangulation the following is true:
Walking along an interior face in clockwise order, the
labels change precisely when moving from a black to a white vertex.
\end{lemma}

\Proof Let $F$ and $F'$ be two faces sharing an edge $e$. Suppose
that $F$ obeys the walking rule. If the clockwise walk in $F$ sees a
change of labels along $e$, then this walk traverses $e$ from the
black to the white vertex, which determines the partition into black
and white for all vertices on $F$ and $F'$. The edge rule implies
that the two labels on the other side of $e$ are the same. This
observation together with the strong face rule for $F'$ yields the
validity of the walking rule for $F'$. The other possibility, when
the clockwise walk in $F$ sees the same label on both ends of $e$,
is similar: The walking rule determines the black/white partition,
the edge rule implies two different labels in $F'$ and the strong
face rule enforces the walking rule for $F'$.

From the definition of labels at the outer face we obtain the
validity of the walking rule for the bounded face $F_0$ which is
incident to the edge containing $s_0$ and having two labels $0$ on
the outer face. Any face $F$ can be connected to $F_0$ with a dual
path avoiding the outer face. The above reasoning allows to transfer
the validity of the walking rule along this path to~$F$. \qed
\bigskip

The following strong edge rule is an immediate consequence of the edge rule together with the
walking rule. Actually, the walking rule also follows from the strong edge rule,
i.e., the two rules are equivalent.

\begin{figure}[htb]
    \begin{center}
        \includegraphics[scale=0.7]{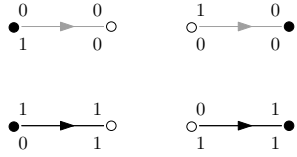}
    \end{center}
    \caption{Edge types complying with the strong edge rule.}
    \label{fig:fourEdgeTypes}
    \end{figure}

\begin{lemma}[Strong edge rule]\label{lem:str-edge}
In a strong labeling of a quadrangulation the following is true: For
each edge, the incident labels coincide at one endpoint and differ
at the other. Moreover if the latter is a white (respectively black)
vertex, the right (respectively left) side of the edge, oriented as
in Observation~\ref{obs:2-coloring_orientation}, has coincident
labels. See Figure~\ref{fig:fourEdgeTypes}.
\end{lemma}

Yet another useful property of strong labelings is given with the
next lemma, whose proof is immediate from the strong edge rule.
Observe that, as in the previous lemma, the rule is
``white--right'', ``black--left''.

\begin{figure}[htb]
    \begin{center}
        \includegraphics[scale=0.9]{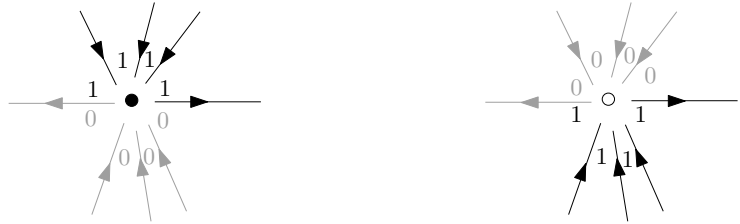}
    \end{center}
    \caption{Illustrating the turning rule.}
    \label{fig:turning-rule}
    \end{figure}

\begin{lemma}[Turning rule]\label{lem:turning}
In a strong labeling of a quadrangulation the following is true:
If~$v$ is a white (respectively black) vertex and~$uv$ an incoming
edge, then the outgoing edge at $v$ with the same color as $uv$ is
the next outgoing edge to the right (respectively left) of $uv$. See
Figure~\ref{fig:turning-rule}.
\end{lemma}

\begin{figure}[htb]
    \begin{center}
        \includegraphics[scale=0.9]{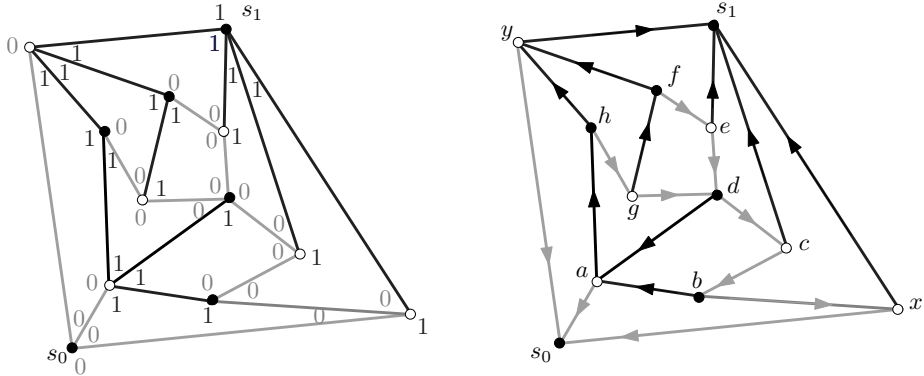}
    \end{center}
    \caption{A strong labeling for a quadrangulation (left) and the induced~2-coloring
    and orientation of the edges (right).}
    \label{fig:quad_trees}
\end{figure}

Figure~\ref{fig:quad_trees} shows a quadrangulation with a strong
labeling. Before further studying strong labelings of quadrangulations
we prove that every quadrangulation has such a labeling.

\begin{theorem}\label{thm:binLabelling}
Every quadrangulation admits a strong labeling.
\end{theorem}

\begin{figure}[htb]
    \begin{center}
        \includegraphics{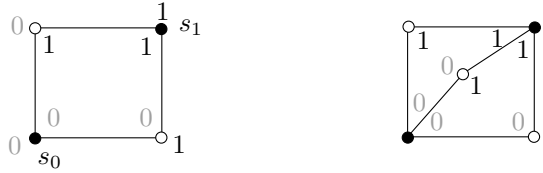}
    \end{center}
    \caption{The basis of the induction and inserting a vertex of degree two.}
    \label{fig:quad_basis_induction}
\end{figure}

\Proof
We use induction on the number of vertices~$n$ of a
quadrangulation~$Q$. If~$n=4$ then a binary labeling exists, as
shown in Figure~\ref{fig:quad_basis_induction} (left). For the
induction step we distinguish two cases.

For the first case, assume that $Q$ contains an interior vertex $v$
of degree two. Removal of $v$ and its two incident edges yields a
quadrangulation $Q^\prime$ which, by induction, admits a binary
labeling. Reinsertion of $v$ and its incident edges into $Q^\prime$
can be done in a unique way such that the rules of strong labelings
are maintained. One of the possible cases is shown in
Figure~\ref{fig:quad_basis_induction} (right).

\begin{figure}[htb]
    \begin{center}
        \includegraphics[scale=0.9]{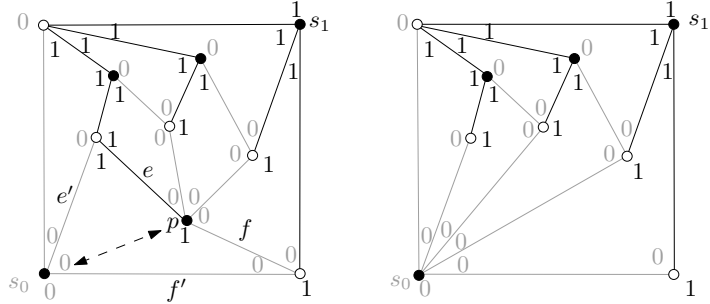}
    \end{center}
    \caption{Contracting a quadrangle to the special vertex $s_0$.}
    \label{fig:quad_contract}
\end{figure}

For the second case, assume that~$Q$ contains no interior vertex of
degree two. We say that a face $q$ incident to $s_0$ is {\emph{contractible}} if it does not contain the
other special vertex $s_1$. The contraction of~$q = \{e^\prime,
e,f,f^\prime\}$, where $\{e^\prime, e,f,f^\prime$\} are the edges of
$q$ in clockwise order starting at~$s_0$ and~$p$ is the vertex
opposite to~$s_0$, identifies $e$ with $e^\prime$, $f$ with
$f^\prime$ and~$p$ with~$s_0$. This can be interpreted as a
continuous movement of $p$ and its incident edges to $s_0$, see
Figure~\ref{fig:quad_contract}.
It is easy to see that if each interior vertex has degree greater than two, then there exists a face incident to the special
vertex~$s_0$ which can be contracted towards~$s_0$.

The contraction of a face $q$ yields a quadrangulation which by induction
admits a binary labeling.
Now, reversing the contraction maintains the binary labeling outside
of the face $q$ and it only remains to label the angles inside $q$.
First, the rule for the special vertex $s_0$ requires that the angle
at this vertex is labeled $0$. The other labels have to
be chosen according to the walking rule. Figure~\ref{fig:quad_contract} shows an example.
The vertex and edge rules
at the boundary of the new face are easily verified.
\qed

From Observation~\ref{obs:strong_quadrangulation} and Theorem~\ref{thm:binLabelling} we deduce
that quadrangulations are precisely the graphs admitting strong labelings.
Nevertheless, for the sake of an easier reading,
in the sequel we will often mention explicitly that the graph considered is a quadrangulation.

\subsection{Schnyder-like properties for strong labelings of quadrangulations}
\label{ssec:schnyder-like_properties}

Consider the coloring and orientation of the edges induced by a
strong labeling of a quadrangulation
(c.f.~Observation~\ref{obs:2-coloring_orientation}). For this
coloring and orientation we obtain results which are in nice
correspondence to those obtained by Schnyder~\cite{Schnyder} for
triangulations and Felsner~\cite{Felsner} for 3-connected plane
graphs. As before, we denote by~$T_i$ the set of oriented edges
colored~$i$ and by $T_i^{-1}$ the set of edges colored $i$ with
reversed orientation.

\begin{lemma} \label{lem:outdegree_1}
Every vertex except~$s_0,s_1$ has
outdegree~1 in each of $T_0$ and $T_1$.
\end{lemma}

\Proof
This follows from the strong edge rule (Lemma~\ref{lem:str-edge}).
\qed
\bigbreak

The next lemma is a special case of
Proposition~\ref{lem:no_directed_cycles-weak}. In particular, it
implies that $T_0$ and $T_1$ are trees.

\begin{lemma}
\label{lem:no_directed_cycles} There is no directed cycle in~$T_0
\cup T_1^{-1}$, nor in $T_1 \cup T_0^{-1}$.
\end{lemma}

\begin{proof}
Although the statement is already settled with
Proposition~\ref{lem:no_directed_cycles-weak} we indicate a second
proof, which is conceptually simpler.
Suppose that there is a cycle in $T_0 \cup T_1^{-1}$. Choose $C$ to
be such a cycle with the least number of faces in its interior.
Claim 1: There is no vertex in the interior of $C$. Otherwise the
black and the gray path leaving the vertex can be used to identify a
cycle with less interior faces. Claim 2: $C$ has no chord. Again
this follows from the minimality assumption.
To complete the proof it can be checked that there is no directed
facial cycle in $T_0 \cup T_1^{-1}$.
\end{proof}

\begin{corollary}\label{cor:trees}
$T_{i}, i \in \{0,1\},$ is a directed tree with sink $s_i$ that
spans all vertices but~$s_{1-i}$.
\end{corollary}

\begin{proof}
Lemma~\ref{lem:no_directed_cycles} implies that~$T_i$ is cycle-free,
hence a forest. Lemma~\ref{lem:outdegree_1} implies that $T_i$ spans
all vertices except~$s_{1-i}$. The same lemma implies that $T_i$ is
directed and has only one sink, $s_i$. This implies the claim.
\end{proof}

For each non-special vertex~$v \notin \{s_0,s_1\}$, we define
the~$i$-path~$P_i(v)$,~$i \in \{0,1\}$, as the directed path
in~$T_i$ from~$v$ to the sink~$s_i$.

\begin{observation}\rm
\label{obs:paths-non-crossing} Two paths of the same color cannot
cross, because every vertex has outdegree 1 in this color. Two paths
of different colors cannot cross, because this would violate the
vertex rule.
\end{observation}

\begin{lemma}\label{lem:chord-free}
The paths $P_i(v), v \notin \{s_0,s_1\},$ are chord-free.
\end{lemma}

\Proof Let $v=v_0,v_1,v_2,\ldots,v_k,s_i$ be the sequence of
vertices of $P_i(v)$. Suppose that $v_iv_j$ with $i+1 < j$ is an
edge of the quadrangulation. The edge is not in the tree $T_i$,
hence, it is of color $1-i$. Lemma~\ref{lem:no_directed_cycles}
implies that the orientation is not from $v_i$ to $v_j$. If $v_jv_i$
lies to the right of $P_i(v)$ we know that $v_j$ is black because of
the turning rule (Lemma~\ref{lem:turning}). The same rule at the
white vertex $v_i$ implies that the outgoing edge at $v_i$ points
into the interior of the cycle $v_i,v_{i+1},\ldots,v_j,v_i$. This
implies a crossing between the paths $P_i(v_i)$ and $P_{1-i}(v_i)$
which contradicts Observation~\ref{obs:paths-non-crossing}. The
other case where $v_jv_i$ lies to the left of $P_i(v)$ is
essentially symmetric. \qed

Because of Observation~\ref{obs:paths-non-crossing}, the paths
$P_{0}(v)$ and $P_{1}(v)$ have $v$ as only common vertex. Therefore
they split the quadrangulation into two regions which we denote
by~$R_0(v)$ and~$R_1(v)$, where~$R_i$ is the region to the right
of~$P_i(v)$ and including both paths.

\begin{lemma}\label{lem:regions}
Let~$u,v$ be distinct interior vertices. For $i \in \{0,1\}$, the
following implications hold:
\begin{itemize}
\item[(i)] $u \in\interior (R_i(v)) \Rightarrow R_i(u) \subset R_i(v).$
\item[(ii)] $u \in P_{i}(v), u\neq v \Rightarrow
\left\{\begin{array}{l}
R_i(u) \subset R_i(v) \text{\quad and \quad} R_{1-i}(v) \subset R_{1-i}(u) \\
\text{or} \\
R_i(v) \subset R_i(u) \text{\quad and \quad} R_{1-i}(u) \subset R_{1-i}(v) \\
\end{array}\right.
$.
\end{itemize}
\end{lemma}

\begin{proof}
If $u \in\interior (R_i(v))$,
Observation~\ref{obs:paths-non-crossing} implies that both paths
$P_0(u), P_1(u) $ and the region they enclose are contained in
$R_i(v)$. If $u \in P_{0}(v)$, $u\neq v$, then $P_{0}(u)\subset P_0(v)$ while
the first edge of $P_1(u)$ points to the interior of either $R_0(v)$
(if~$u$ is black) or to the interior of $R_1(v)$ (if~$u$ is white),
because of the turning rule. In the first case we obtain $R_0(u)
\subset R_0(v)$ and $R_1(u) \supset R_1(v)$, in the second case we
obtain the reversed inclusions. Similar arguments work if $u \in
P_{1}(v)$, $u\neq v$.
\end{proof}

\subsection{Alternating embedding of quadrangulations on 2-books}
\label{ssec:2-book_embedding}

Mimicking the obtention of straight-line embeddings of triangulations
on small grids via Schnyder labelings, Lemma~\ref{lem:regions} allows us to
obtain $2$-book embeddings of quadrangulations such that each page
contains an alternating tree.

For each non-special vertex~$v$, let us define~$f_i(v)$ as the number of
faces contained in~$R_i(v)$. For the two special vertices~$s_0,
s_1$, we set $f_0(s_0) = f_1(s_1) = -1$ and $f_1(s_0)= f_0(s_1) =
n-2$, where~$n$ is the number of vertices. As shown in
Lemma~\ref{lem:regions} there is an inclusion between the
$i$-regions of any two vertices, therefore, the following holds:

\begin{proposition}
For any two vertices $u \neq v$, we have $f_i(u) \neq
f_i(v)$. Equivalently, all possible values of $f_i$ from $0$ to $n-3$ occur.
\end{proposition}

All the points $(f_0(v), f_1(v))$ lie equally spaced on the
line~$f_0+f_1=f$, where~$f$ is the total number of bounded faces
(which equals~$n-3$ by Euler's formula). For the sake of
convenience, we can choose a reference system in which this line is
the horizontal axis and \emph{the~$f_1$-values increase from left to
right}. Given this as spine of the book, we draw the edges of each tree~$T_i$ on
one side. As a convention, we will draw~$T_0$ gray and above the
line, and~$T_1$ black and below. In Theorem~\ref{thm:non-crossing}
we prove that the trees are non-crossing, and hence we get
a~$2$-book embedding for the quadrangulation~$Q$ such that each page
contains a tree.

\begin{figure}[htb]
\centering
\includegraphics[scale=0.9]{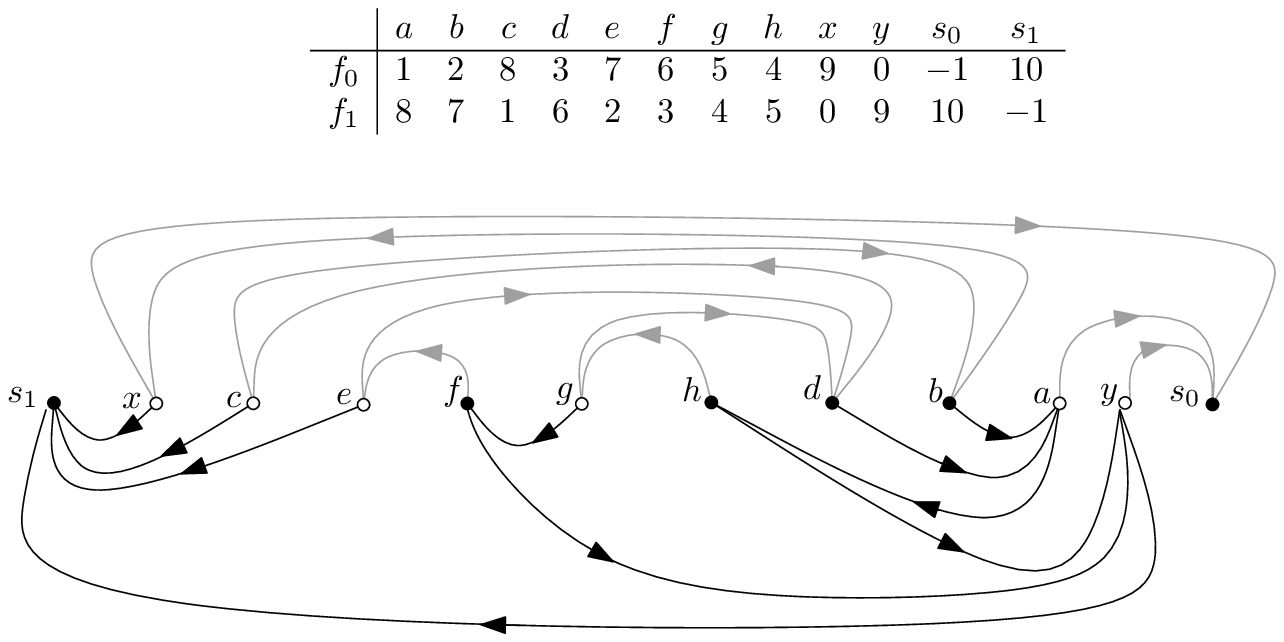}
    \caption{Embedding on a~$2$-book the quadrangulation in Figure~\ref{fig:quad_trees}.}
    \label{fig:quadrangulation_labelling_lineembedding}
\end{figure}

In Theorem~\ref{thm:non-crossing} we additionally prove that both
trees are {\emph{alternating}}, meaning that the tree contains no
two edges~$v_iv_j$ and~$v_jv_k$ for~$i<j<k$ (where $v_1,\dots,v_n$ denotes the vertices in the order they are encountered along the line). This is equivalent to
saying that either all neighbors of $v_j$ have indices bigger than
$j$ or they all have indices smaller than $j$.
Figure~\ref{fig:quadrangulation_labelling_lineembedding} shows an
example for the book embedding.

Non-crossing alternating trees are counted by the Catalan numbers.
They came up in research about pseudo-triangulations, where
they have been identified as one-dimensional analogs to
pseudo-triangulations~\cite{RSS}. In that paper it has been
shown that the ``flip graph'' on alternating trees is
the $1$-skeleton of the associahedron.

\begin{theorem}\label{thm:non-crossing}
Let the vertices of a quadrangulation be placed on a line by the face-counting process, with the
trees $T_0$ and $T_1$ placed on each side of the line. Then $T_0$ and $T_1$ are non-crossing and alternating.
\end{theorem}

\begin{figure}[htb]
\centering
\includegraphics[scale=0.5]{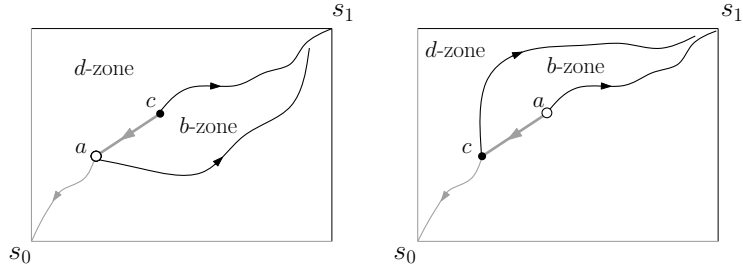}
\caption{Possible configurations according to the relations
in~(\ref{eq:f_0}).} \label{fig:PlanarityCases}
\end{figure}

\Proof
We will prove that the gray tree~$T_0$ cannot have crossings.
Let us
suppose that there is a crossing in~$T_0$, i.e., four points
$a,b,c,d$ with

\begin{equation}\label{eq:f_0}
f_0(a)>f_0(b)>f_0(c)>f_0(d)
\end{equation}
and edges~$ac,bd$. We focus on the edge~$ac$. The two possible
configurations are shown in Figure~\ref{fig:PlanarityCases}. Any
other situation would violate either the relations in~(\ref{eq:f_0})
or the vertex rule~(G1).

Furthermore, from~$f_0(c)>f_0(d)$ and Lemma~\ref{lem:regions}, we
know that~$d\not\in\interior(R_1(c))$ and analogously
that~$b\not\in\interior(R_0(c))$. This gives us the feasible zones
for points~$b$ and~$d$, denoted as~$b$- and~$d$-zones in
Figure~\ref{fig:PlanarityCases}. Note that in both cases shown in the
figure the path $P_1(c)$ separates the two zones, hence, the existence of a gray
edge $bd$ implies that $b$ or $d$ is on this path.
From Lemma~\ref{lem:chord-free} we know that at most one of them is
on the path.

\begin{itemize}
\item {If $b\in P_1(c)$, $d\not\in P_1(c)$:}
In this case the edge~$bd$ is
to the left of~$P_1(c)$ and hence the same holds for all the gray
edges incident to~$b$. Therefore, $R_0(b) \subset R_0(c)$ which
implies $f_0(b) < f_0(c)$, a contradiction.

\item {If $b\not\in P_1(c)$, $d\in P_1(c)$:} In this case the edge~$bd$ is
to the right of~$P_1(c)$ and hence the same holds for all the gray
edges incident to~$d$. This leads to the contradiction
$f_0(d)>f_0(c)$.

A similar analysis shows that the black tree~$T_1$ has no crossings.
\end{itemize}

We now show that for our choice of coordinates,~$T_0$ and~$T_1$ are alternating.
We focus on the black tree~$T_1$, and the case for the gray
tree~$T_0$ is analogous. Vertices incident to the exterior face have all their
neighbors on one side, hence they are alternating.

Let us consider an interior black vertex~$v_j$.
The successor $v_s$ of $v_j$ on the black path $P_1(v_j)$ is a white
vertex. From the turning rule (Lemma~\ref{lem:turning}) it follows that
the gray outgoing edge of $v_s$ is to the left of $P_1(v_j)$, i.e., it
points into $R_0(v_j)$ which implies $f_0(v_s) < f_0(v_j)$,
equivalently $f_1(v_s) > f_1(v_j)$ and hence $s > j$.

Now consider a black edge $v_pv_j$ which is incoming at $v_j$. This edge
belongs to the black path $P_1(v_p)$. The fact that $v_j$ is black and
the turning rule implies that the gray outgoing edge of $v_j$ points
into $R_1(v_p)$ which implies  $f_1(v_p) > f_1(v_j)$ and hence $p > j$.

The case where $v_j$ is a white vertex is similar, in that case all
neighbors in the black tree have indices smaller than $j$.
\qed

\subsection{Strong labelings, separating decompositions and 2-orientations
}
\label{sec:separating_decomp}

The following definition was essentially (with reversed orientations)
proposed by de~Fraysseix and Ossona~de~Mendez~\cite{Fraysseix}.

\begin{definition}\rm
Let~$Q$ be a quadrangulation, with vertices of the bipartition
properly bicolored as black and white. Let $s_0$ and $s_1$ be
nonadjacent vertices at the outer face. A
\emph{separating decomposition} of~$Q$ is a partition of the edges
into two directed trees~$T_0,T_1$ with sinks~$s_0,s_1$,
such that the incident edges at each vertex but $s_0$ and $s_1$ are
gathered as follows, in clockwise order for black vertices and
counterclockwise order for white vertices:
\Bitem{The incoming edges (if any) from $T_0$,}
\Bitem{The outgoing edge from $T_0$,}
\Bitem{The incoming edges (if any) from $T_1$,}
\Bitem{The outgoing edge from $T_1$.}
\end{definition}

\ni
Note that the above condition about the orientations of edges at a vertex
is exactly the turning rule (see Figure~\ref{fig:turning-rule}).

\begin{theorem}\label{thm:bij sd-gl}
Separating decompositions and strong labelings of a quadrangulation
are in bijection.
\end{theorem}

\Proof
Let $Q$ be a quadrangulation with a  distinguished
vertex $s_0$ on the outer face. A strong labeling of $Q$
induces a coloring and orientation of the edges. By
Corollary~\ref{cor:trees} this
yields a partition into trees $T_0$ and $T_1$ rooted at $s_0$ and
the opposite vertex of the outer face $s_1$. The
coloring and orientation of the edges obeys the turning rule
(Lemma~\ref{lem:turning}). This rule is precisely the
condition required for a separating decomposition.

Conversely, let a separating decomposition be given.
Given a directed edge $uv$ color both angles incident to $uv$ at $v$
with the color of the edge. The separation property implies that
angles with two incident incoming edges get the same label from both edges.
Angles which are unlabeled at this point are labeled according to
the strong edge rule (see Figure~\ref{fig:fourEdgeTypes}).
It is obvious that the vertex conditions (G0) and (G1) hold
for this labeling. All edges conform to the strong edge rule
and hence the edge rule (G2). The strong edge rule also
implies the walking rule (Lemma~\ref{lem:walking})
which in turn implies the strong face rule (G3$^+$).
Together this shows that the implied labeling of angles
is a strong labeling.
\qed
\bigskip

To enhance the picture we quote the following theorem
from~\cite{Fraysseix}. In the statement `quadrangulation' is again to
be understood as a quadrangulation together with a distinguished
vertex $s_0$ on the outer face.

\begin{theorem}[De~Fraysseix and Ossona~de~Mendez] \label{thm:bij sd-2or}
Separating decompositions and 2-orientations of a quadrangulation
are in bijection.
\end{theorem}

\begin{corollary}\label{cor:2or-strong}
There is a bijection between 2-orientations and strong labelings of a quadrangulation.
\end{corollary}

\begin{figure}[htb]
    \begin{center}
        \includegraphics{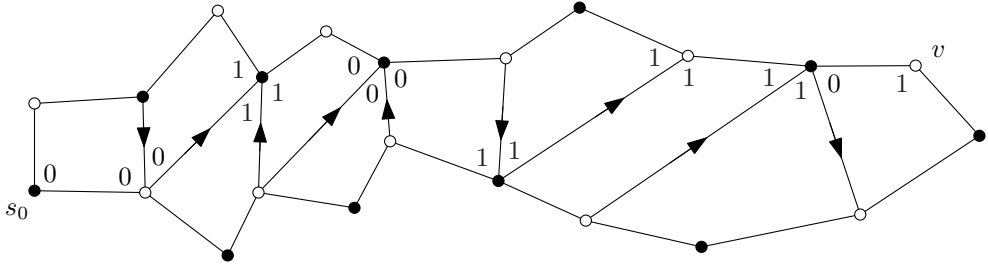}
    \end{center}
    \caption{The label at $v$ is deduced from the label at $s_0$ using the oriented edges and the walking rule.}
    \label{fig:transfer-labeling}
\end{figure}

\Proof
Theorems~\ref{thm:bij sd-gl} and~\ref{thm:bij sd-2or} imply
that : We already know how a strong labeling, actually even a
weak labeling, induces a 2-orientation. For the reverse mapping, suppose
a 2-orientation of $Q$ with $s_0$ is given and we want to find the
corresponding strong labeling. Note that by the walking rule
(Lemma~\ref{lem:walking}) we know all the labels of a face if we
know just one. If the labels on one side of an oriented edge are
given, then we can copy the label at the tip of the edge to the other side
and deduce all labels in the face of that side by the walking rule.
This allows to infer all angle labels of an 2-oriented quadrangulation
just from the labels at $s_0$.
An example is given in Figure~\ref{fig:transfer-labeling}. In order to complete
 the proof it just remains to check that the above
method will not yield a conflicting assignment of labels to an angle
and that the resulting labeling has the properties required for a
strong labeling.
\qed
\bigskip

\subsection{Flips on strong labelings of quadrangulations}
\label{subsec:flips}

\begin{figure}[htb]
    \begin{center}
        \includegraphics[scale=0.8]{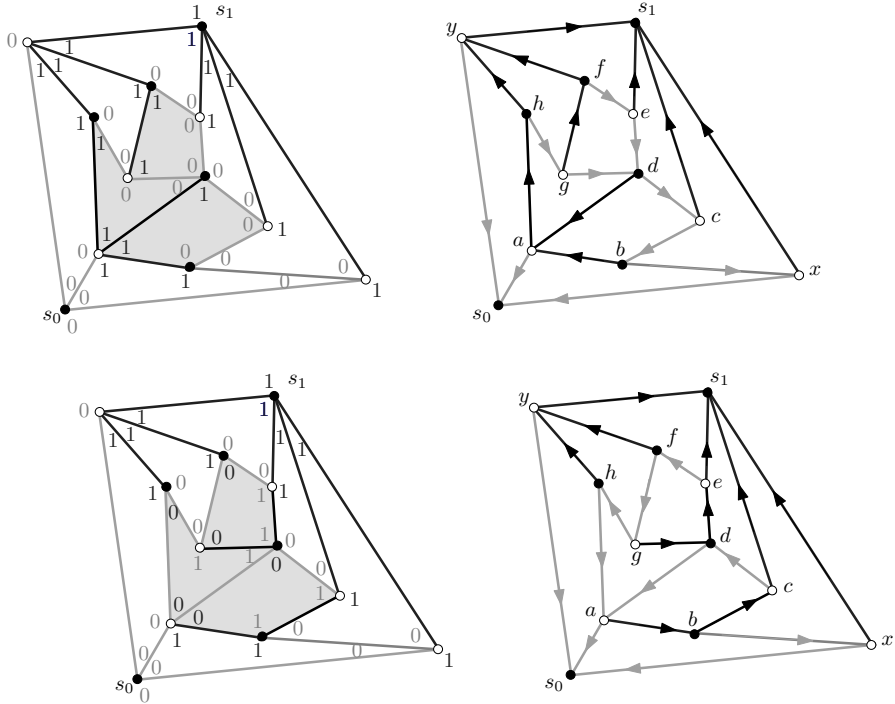}
    \end{center}
    \caption{A flip for a strong labeling of a quadrangulation.}
    \label{fig:quadr_labeling_flip}
\end{figure}

Given a graph~$G$ with a 2-orientation, one can reverse any directed cycle and obtain
another 2-orientation. Such a local modification of an object
in more general terms is often called \emph{flip}.
In Lemma~\ref{lem:flip} below we show that for a strong labeling of a quadrangulation, such a flip
means that we invert all labels inside the cycle to obtain
another strong labeling of the same quadrangulation, see
Figure~\ref{fig:quadr_labeling_flip}.

\begin{lemma}\label{lem:flip}
Given a quadrangulation, the reversal of a directed cycle $C$ in a 2-orientation conforms
with the complementation of all the labels inside $C$
in the corresponding strong labeling.
\end{lemma}

\begin{proof}
Recall the description of how to compute the strong labeling
corresponding to a 2-orientation (Figure~\ref{fig:transfer-labeling}).
From the edge rule it follows that whenever an edge is used to infer
the labels in a face the reoriented edge will imply the complementary
labels in that face. The computation of labels starts at $s_0$ which is
at the outer face, hence, outside of $C$.  The deduction of the label
of an angle outside of $C$ always uses an even number of edges of $C$,
hence, the label is not changed by reorienting $C$. But the
computation of the label of an angle inside $C$ uses an odd number of
edges of $C$ and the label is complemented.
\end{proof}

Schnyder woods of triangulations are in bijection to 3-orientations.
In this context Brehm~\cite{Brehm} has investigated the reversal of
directed cycles (flip) for 3-orientations. He proved that the set of
3-orientations forms a distributive lattice. More generally Ossona
de Mendez~\cite{mendez} and Felsner~\cite{Felsner_latticeStruc}
found lattice structures on the set of $\alpha$-orientations of a
planar graph. A particular instance of the general theorem is that
the set of all 2-orientations of a quadrangulation can be enhanced
with an ordering which is a distributive lattice. The order relation
is the transitive closure generated by $X < X_C$ whenever $X$ is a
2-orientation which has a simple directed cycle $C$ which runs
clockwise around its interior and $X_C$ is obtained by reverting $C$
in $X$. The flip structure on 2-orientations of quadrangulations was
also investigated by Nakamoto and Watanabe~\cite{NW}. A simple
consequence of the distributive lattice structure is:

\begin{corollary}
The flip graph of strong labelings is connected.
\end{corollary}

\section{Generalized strong labelings}\label{sec:genstronglab}

In Section~\ref{sec:binary_labelings} we introduced weak labelings, which can only exist for plane graphs with~$n$ vertices and~$2n-4$ edges. Strengthening the face rule in Section~\ref{ssec:strong-quadrangulations} resulted into
quadrangulations, i.e., maximal bipartite plane graphs, as the only graphs that admit strong labelings, i.e., conditions (G0), (G1), (G2) and~(G3$^+$).
In this section we modify the edge rule~(G2) in order to have similar labelings for a larger class of bipartite plane graphs. The following is inspired by the generalization~\cite{Felsner} of Schnyder woods for 3-connected plane graphs.
We will always assume that one color class of a bipartite graph has been selected to be the white class,
the other one is the black class.

\begin{definition}\rm
A {\emph{generalized strong labeling}} for a bipartite plane graph is a
mapping from its angles to the set $\{0,1\}$ which satisfies (G0),
(G1), (G2$^+$) and (G3$^+$), where

\Item(G2$^+$) {\bf{Extended edge rule:}} For each edge, the incident labels
     form one of the six patterns shown in Figure~\ref{fig:sixEdgeTypes}.

    \begin{figure}[htb]
    \begin{center}
        \includegraphics[scale=0.8]{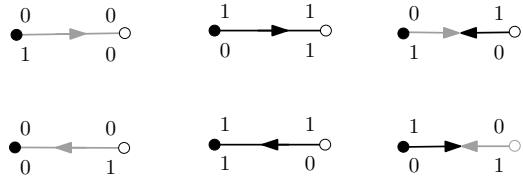}
    \end{center}
    \caption{Extended edge rule.}
    \label{fig:sixEdgeTypes}
    \end{figure}

\end{definition}

Figure~\ref{fig:strong-ex} shows several examples of generalized strong labelings.
This definition is justified by the following result:

\begin{lemma}
A generalized strong labeling of a quadrangulation
has only edges of the four types of the edge rule~(G2), which verify
Lemma~\ref{lem:str-edge}, i.e., it has no bidirected edge.
\end{lemma}

\Proof
A quadrangulation on $n$ vertices has $n-2$ faces and $2n-4$ edges.
Every face requires two color changes in its face walk. Every edge contributes
at least one color change to a face walk.
\qed
\bigskip

It should be noted that the degree of the special vertices $s_0$ and $s_1$
in a graph with such a labeling can be one, e.g., the 2-path $s_0$---$v$---$s_1$
admits a generalized strong labeling.

\begin{figure}[htb]
    \begin{center}
        \includegraphics[scale=0.6]{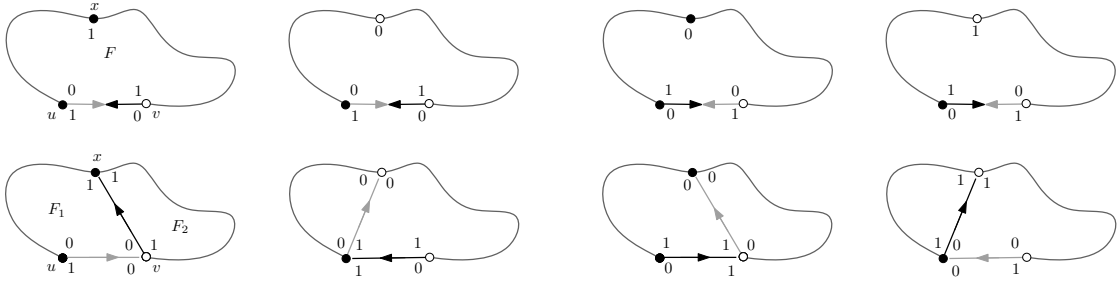}
    \end{center}
    \caption{Split and merge for generalized strong labelings.}
    \label{fig:merge+split}
\end{figure}

Bonichon et al.~\cite{BoFeMo} have introduced operations on Schnyder
woods which they call \emph{merge} and \emph{split}. A split takes a
bidirected edge and opens it up into two unidirected edges. A merge
is the inverse operation; it takes an angle with two unidirected
edges, one of them incoming the other outgoing, and turns the
outgoing edge into the incoming thus making it bidirected. We define
similar operations for generalized strong labelings.
Figure~\ref{fig:merge+split} shows the four possible instances for
split and merge.  A split is done by replacing a situation from the
upper row by the situation below. A merge, conversely, replaces a
situation in the lower row by the one above it.

\begin{lemma}
If $G$ is a graph with a generalized strong labeling $B$ and a labeling $B'$ of
$G'$ is obtained from $(G,B)$ by a split or merge, then $B'$ is a
generalized strong labeling of $G'$.
\end{lemma}

\Proof Assuming that $B$ obeys the vertex rules (G0) and (G1) these
rules are easily seen to hold for $B'$ as well. All edges in the
figure are legal in the sense of (G2$^+$). The least trivial thing
is to verify (G3$^+$) for the split. Let us concentrate on the split
of the first column where we have given names to the objects. The
two black vertices $u$ and $x$ have different labels inside $F$.
Hence, when walking clockwise from the edge $vu$ towards $x$ we have
to pass at exactly one of the two edges which have identical labels
at both ends inside $F$, by rule (G3$^+$). Before reaching this edge
we always see a 0 at black vertices and a 1 at white vertices. From
Figure~\ref{fig:sixEdgeTypes} we find that a clockwise traversal of
an edge with identical labels always goes from the white to the
black vertex. Hence, the edge we meet has two labels 1. This is what
we need to show that (G3$^+$) holds for $F_1$. Similar arguments
show that (G3$^+$) holds for $F_2$ and indeed that it holds for the
two new faces after each of the four possible splits. \qed

\begin{lemma}\label{lem:split}
Let $G$ be a graph with a generalized strong labeling. If $G$ is not a
quadrangulation then there is an edge which is feasible for a split.
\end{lemma}

\Proof If $G$ is not a quadrangulation then it has more edges than
twice the number of faces. Therefore there is a bidirected edge
$uv$. Let $u$ be black and $v$ be white and consider the face $F$
whose clockwise traversal sees $e$ as the edge from $v$ to $u$. We
assume that the label of $v$ in $F$ is 1. From the proof of the
previous lemma we deduce that clockwise from $vu$ we reach the edge
with labels $1,1$ and that the second vertex $x$ of this edge is
black. This shows that a split of the edge $uv$ towards $x$ is possible.
The case in which the label of~$v$ in~$F$ is~$0$ works analogously.

A special case occurs if the face $F$ is the outer face. To handle
this case think of $G$ as being embedded on the sphere and note that
the special conditions of (G3$^+$) for the outer face impose the same
structure we have noted for the other faces. Hence splits are possible
but special care must be put into the choice of the vertex $x$ towards
which an edge is split, a careless choice could split the outer face
such that there is no face containing both $s_0$ and~$s_1$.
\qed

\begin{corollary}\label{cor:split2Q}
If $G$ is a graph with a generalized strong labeling then there is a sequence of
edge splits which lead to a quadrangulation with a strong labeling.
\end{corollary}

This corollary has quite strong consequences, e.g., we may observe that
the turning rule (Lemma~\ref{lem:turning}) is invariant under splitting and
merging. Hence the turning rule holds for graphs with a generalized strong labeling.

Given a graph $G$ with a generalized strong labeling, let $T_0$ be the set
of oriented gray edges and let $T_1$ be the set of oriented black edges.
Again $T_i^{-1}$ is the set of edges of $T_i$ with reversed orientation.

\begin{lemma} \label{lem:trees}
$T_0 \cup T_1^{-1}$ and $T_1 \cup T_0^{-1}$ are acyclic.
Moreover, $T_{i}, i \in \{0,1\},$ is a directed tree with sink $s_i$ that
spans all vertices but~$s_{1-i}$.
\end{lemma}

\Proof Use edge splits to get from $G$ to a quadrangulation $Q$. The
acyclicity of $T_0 \cup T_1^{-1}$ where $T_i$ are the edge sets
defined by the orientation of $Q$ was shown in
Lemma~\ref{lem:no_directed_cycles}. Note that since a merge has
precisely the effect of deleting an edge from $T_0 \cup T_1^{-1}$,
this cannot introduce cycles.

The statement about the trees again follows from the acyclicity of
$T_i$ and the fact that every non-special vertex has outdegree one in
$T_i$.
\qed
\medskip

The lemma implies that, again, we can define the~$i$-path~$P_i(v)$,
of a vertex $v$ as the directed path
in~$T_i$ from~$v$ to the sink~$s_i$. These paths allow, in turn,
the definition of the regions $R_i(v)$ of a vertex.
Consider the numbers $f_0(v)$ counting the number of faces in $R_0(v)$,
i.e., the region to the right of~$P_0(v)$. These numbers again
obey a nice {\em alternation property}, namely, if $xy$ is an edge
of color 0 with black end~$x$ and white end~$y$, then $f_0(x) \leq f_0(y)$.
If the color of $x,y$ is 1 and $x$ is black and $y$ white, then $f_0(x) \geq f_0(y)$.
However, we lose the property that the numbering $f_0$ yields a 2-book embedding;
this is due to the fact that a 0-path $P_0(u)$ and a 1-path $P_1(v)$ can cross
by using the two directions of a bidirected edge.

\subsection{Distributive lattice and flips for generalized strong labelings}
\label{ssec:lattice_flips_genstronglab}

In Section~\ref{ssec:strong-quadrangulations} we showed that strong
labelings for quadrangulations are in bijection to 2-orientations.
This allowed us to identify a flip operation on strong labelings
which generates a distributive lattice on the set of all strong
labelings. The following construction allows to prove equivalent
results in the case of generalized strong labelings.

\begin{figure}[htb]
    \begin{center}
        \includegraphics[scale=0.9]{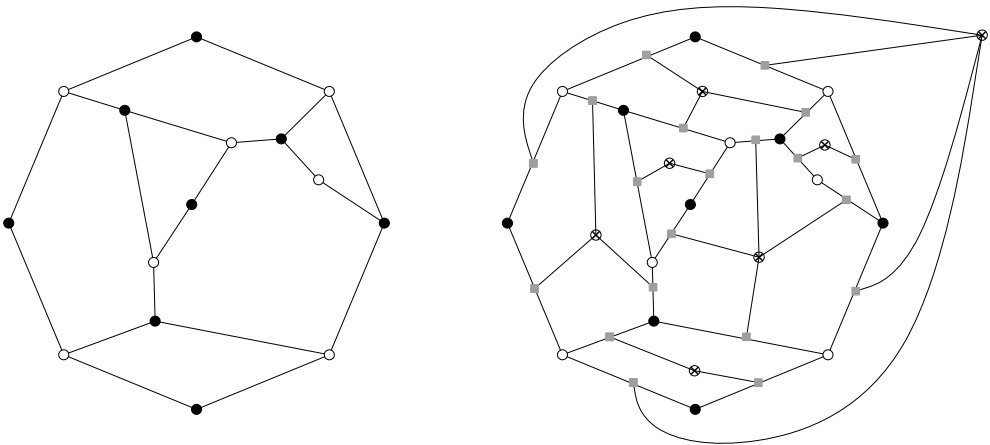}
    \end{center}
    \caption{A graph $G$ (left) and the corresponding $S_G$ (right).}
    \label{fig:half-completion}
\end{figure}

The orientation induced by a generalized strong labeling on $G$ has the somewhat
strange property that it may contain bidirected edges. We encode
this orientation by a ``regular'' orientation of a bigger graph: Let
$G$ be a connected bipartite plane graph with distinguished color
classes black and white and two special vertices $s_0$ and $s_1$ on
the outer face. Define a graph $S_G$ as follows: As vertices of
$S_G$ take the union of the vertices, edges and faces of $G$. Every
edge-vertex has degree three and is connected to the two endpoints
and to the face on its right when traversed from white to black.
Figure~\ref{fig:half-completion} shows an example. The construction
somewhat resembles the completion of a plane graph as
used in the proof of Proposition~\ref{prop:weak_lab+orient}.
Similar constructions have been considered in the context of Schnyder wooods,
see e.g.~\cite{Felsner_latticeStruc}.

\begin{proposition}
Generalized strong labelings of $G$ are in bijection with orientations of $S_G$
which have the following outdegrees
$$
{\sf outdeg}(x) =
       \begin{cases}
             0 & \text{if $x\in\{s_0,s_1\}$,}\\
             1 & \text{if $x$ is an edge-vertex,}\\
             2 & \text{otherwise.}
       \end{cases}
$$
\end{proposition}

\begin{figure}[htb]
    \begin{center}
        \includegraphics{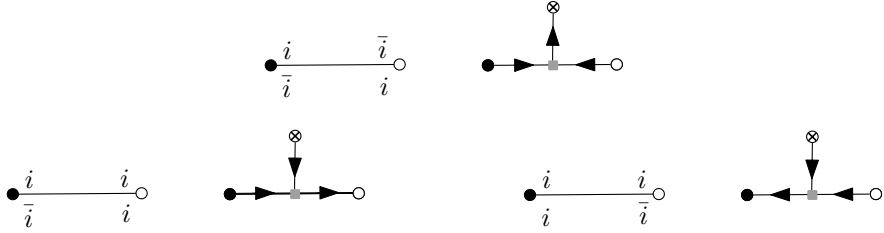}
    \end{center}
    \caption{Translating orientations from $G$ to $S_G$ and back. $\bar{i}$ denotes the label $1-i$.}
    \label{fig:orientations}
\end{figure}

\Proof Figure~\ref{fig:orientations} shows how to translate from a
generalized strong labeling of $G$ to an orientation of $S_G$. 
There is a clear correspondence between the rules (G0) and (G1) and
the prescribed outdegrees of original vertices. The extended edge rule and
outdegree 1 for edge-vertices are both assumed for the
translation. The face rule (G3$^+$) corresponds to
outdegree 2 for face-vertices. Note that this also holds for the
outer face, the two edges on the outer face which should have
repeated labels to satisfy (G3$^+$) connect to
the vertices $s_0$ and $s_1$ which have prescribed outdegree $0$.
Therefore, these two edge-vertices receive the two outgoing edges
of the vertex of the outer face.
\qed
\bigskip

The orientations of $S_G$ described in the proposition are
$\alpha$-orientations in the sense of~\cite{Felsner_latticeStruc}.
Hence, the set of all generalized strong labelings of $G$ can be ordered as a
distributive lattice. In particular the generalized strong labelings are again
flip-connected, where a flip is defined as the complementation of
all labels inside a cycle $C$ which is directed in the corresponding
orientation of $S_G$.

Recall that for given~$G$ and a fixed function~$\alpha$, the existence
of an $\alpha$-orientation can be decided in polynomial time.
Together with the proposition, this yields:

\begin{theorem}
Plane graphs admitting a generalized strong
labeling can be recognized in polynomial time.
\end{theorem}

\subsection{Graphs admitting a generalized strong labeling}
\label{ssec:admit_genstronglab}

So far we have shown that generalized strong labelings have a nice structure.
The correspondence with orientations of $S_G$ yields polynomial
time recognition and an implicit characterization via the
criterion for $\alpha$-orientations given on page~\pageref{demand_crit}. In this
subsection we provide an explicit characterization.
Together with the classical algorithm of Hopcroft and Tarjan~\cite{HoTa}
for finding the triconnected components this allows a linear time
recognition of the class.

\begin{figure}[htb]
    \begin{center}
        \includegraphics{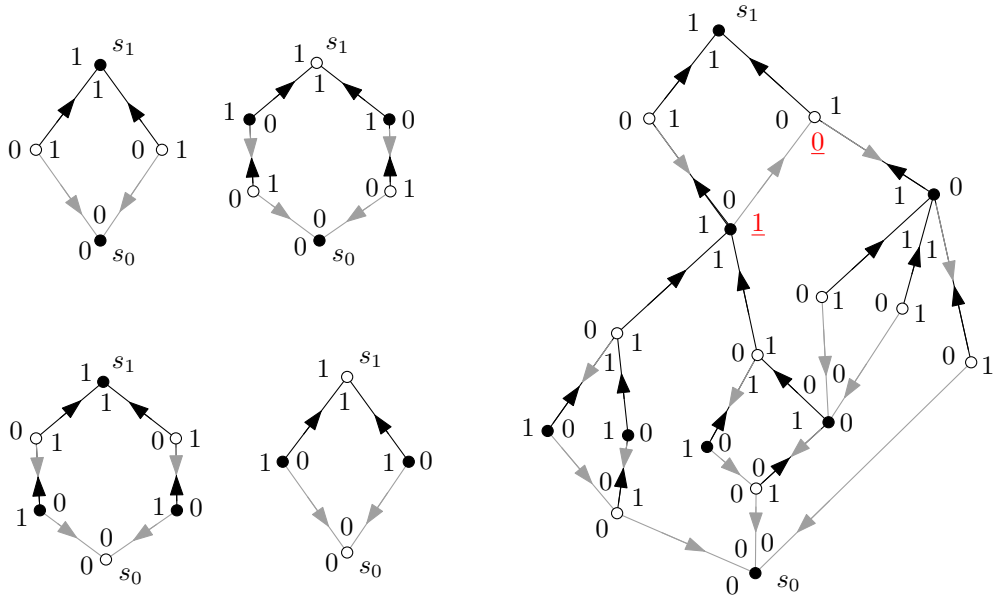}
    \end{center}
    \caption{Examples of generalized strong labelings.}
    \label{fig:strong-ex}
\end{figure}

To introduce into the topic we have two figures. Figure~\ref{fig:strong-ex}
shows some examples of graphs with generalized strong labelings. The four examples on the
left illustrate how the colors of the special vertices influence the labeling
along the outer face. The generalized strong labelings in these cases are unique. The generalized strong
labeling of the larger graph on the right is not unique, e.g., exchanging the
two underlined labels leads to another generalized strong labeling.
\begin{figure}[htb]
    \begin{center}
        \includegraphics{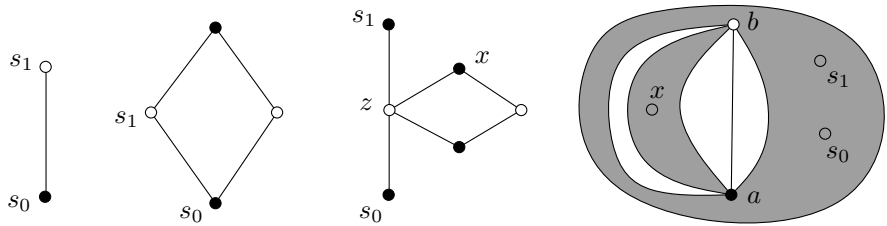}
    \end{center}
    \caption{Some graphs which do not admit a generalized strong labeling.}
    \label{fig:bad-ex}
\end{figure}
Figure~\ref{fig:bad-ex} shows some graphs which do not admit generalized strong
labelings for different reasons. The first two examples fail to
admit a generalized strong labeling simply because their two special vertices
are adjacent. Rule (G0) would force the connecting edge to have two
identical labels on both ends, which is infeasible by the edge rule.
In the middle example there is a cut vertex between $x$ and the two
special vertices. The two paths $P_0(x)$ and $P_1(x)$ would both
contain $z$ which forces a cycle in $T_0 \cup T_1^{-1}$, which is
impossible by Lemma~\ref{lem:trees}.

An undirected graph with special vertices $s_0$ and $s_1$ is called
\emph{weakly 2-connected} if it is 2-connected or adding an edge
$s_0s_1$ makes it 2-connected. This is equivalent to saying that every
vertex $x$ has a pair of vertex-disjoint paths one leading to $s_0$ and the
other to $s_1$. From the above it follows that being weakly
2-connected is a necessary condition for admitting a generalized strong labeling.

Now consider the sketch on the right of Figure~\ref{fig:bad-ex}. It
illustrates the following situation: There is an edge $ab$, vertex
$a$ is black and vertex $b$ white.  Removing $a$ and $b$ we
disconnect a component $C$ with $x\in C$ from the special vertices
$s_0$ and $s_1$. Moreover, component $C$ is to the left of $ab$. If
a graph contains such an edge we say that it \emph{contains a block
with a right chord}. Suppose that a graph containing a block with a
right chord admits a generalized strong labeling. Disjointness forces the two
paths $P_0(x)$ and $P_1(x)$ to leave the component $C$ through
vertices $a$ and $b$. From Lemma~\ref{lem:trees} it can be concluded that
there is no edge oriented from $a$ or $b$ into $C$.  Now consider
the orientation of the edge $ab$, if it is directed from $b$ to $a$,
then the turning rule for white vertices makes the path $P_i(x)$
leaving at $b$ continue through $a$ where the two paths meet,
contradiction. If the direction of $ab$ is from $a$ to $b$, then it
is the turning rule for the black vertex $a$ which leads to the same
kind of contradiction.

With the three cases we have identified all the obstructions against
admitting a generalized strong labeling:

\begin{theorem}
Let $G$ be a bipartite plane graph with color classes black and
white and two special vertices $s_0$, $s_1$ on the outer face. $G$
admits a generalized strong labeling if and only if the following conditions
are satisfied. \Item(1)  $s_0$ and $s_1$ are nonadjacent, \Item(2)
$G$ is weakly 2-connected, \Item(3)  $G$ contains no block with a
right chord.
\end{theorem}

\begin{figure}[htb]
    \begin{center}
        \includegraphics[width=\textwidth]{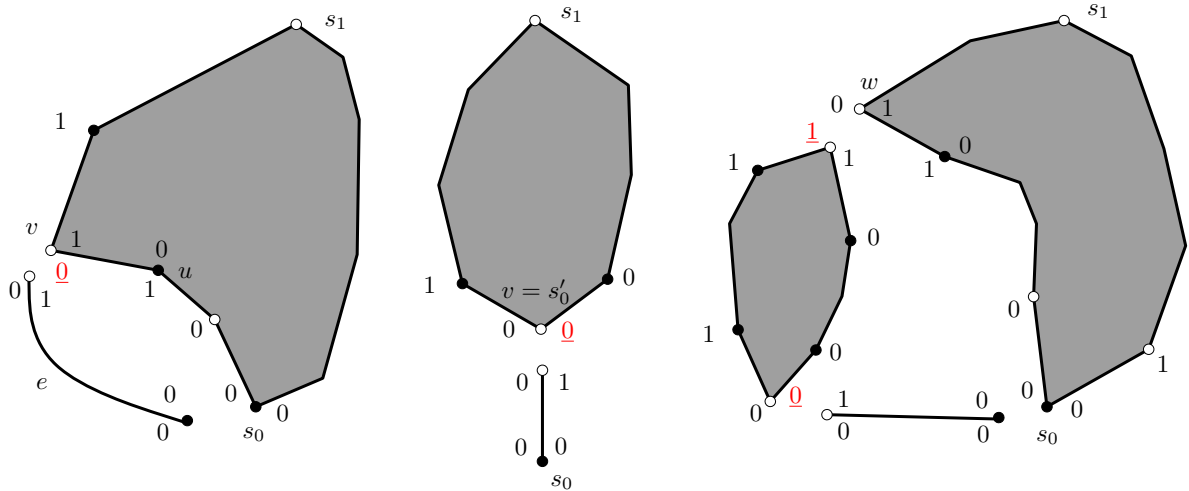}
    \end{center}
    \caption{Constructing the generalized strong labeling in the inductive proof. Underlined
    labels are inverted in the labeling of $G$.}
    \label{fig:charact-proof}
\end{figure}

\Proof The ``only if'' part comes from the above discussion. The
proof for the ``if'' part is by induction on the number of edges.
Let $G$ be a graph satisfying the conditions. We concentrate on the
case where $s_0$ is a black vertex, and the other case is similar.
Let $e=s_0v$ be the first edge in clockwise order which is incident to~$s_0$
and belongs to the boundary of the outer face
(in Figure~\ref{fig:strong-ex} it is the leftmost edge at $s_0$). Rule~(G2$^+$)
implies that~$e$ has the duplicate label $0$ on the outer
face. Now, remove $e$ from $G$ and let $G'$ be the resulting graph. There are several
cases, Figure~\ref{fig:charact-proof} shows how to deal with them.

The first case is that $G'$ satisfies the conditions and we can by
induction assume a generalized strong labeling for $G'$. Consider the edge $uv$
on the boundary of the outer face of $G'$ which is interior in $G$.
In the labeling of $G'$ on the outer face the black vertex $u$ has
label $1$ and the white vertex $v$ has label $0$. The extended edge
rule (G2$^+$) implies that the labels on the opposite side of this
edge are inverse, $0$ at $u$ and $1$ at $v$. Therefore, it is
consistent with edge and vertex rules to label the angle between $e$
and $uv$ with $1$ and the outer angle of $e$ at $v$ with $0$. This
yields a generalized strong labeling of $G$.

If $G'$ does not satisfy the conditions then, necessarily, it is
condition (2) which fails. If $G'$ is not connected it has $s_0$ as
an isolated vertex. Choose $v$ as the special vertex $s'_0$ for the
component of $G'$ which contains $s_1$. If this component admits a
generalized strong labeling we can extend this to the full graph. Otherwise,
condition (1) is not satisfied. Hence either the component is just
the single edge $s_0's_1$ or this edge is a left chord to a block
which satisfies all three conditions. In both cases it is easy to
get to a generalized strong labeling of $G$.

If $G'$ is connected but fails to satisfy (2), then it has a cut
vertex. Let $w$ be the cut vertex such that one of the components is
weakly 2-connected between $v$ and $w$ and the other is weakly
2-connected between $s_0$ and $s_1$. The first of these components
is either a single edge or it satisfies the conditions. The second
component also satisfies the conditions. By induction both
components have generalized strong labelings. Again it is straightforward to
define a generalized strong labeling based on the generalized strong labelings of the
components. The right part of Figure~\ref{fig:charact-proof} shows
the case where $w$ is a white vertex. \qed

\section{A binary labeling for plane Laman graphs}
\label{sec:A binary labeling for plane Laman graphs}

Although a deep study is left for further work, let us point out in this section that weak labelings can be extended as well to \emph{Laman graphs}, those with $n$ vertices and $2n-3$ edges such
that any induced subgraph on $k$ vertices has at most $2k-3$ edges.
These graphs arise in the context of rigidity theory~\cite{Laman} and they are strongly connected to pseudo-triangulations~\cite{osss-cpt-2007,Streinu}, i.e., plane straight-line drawings such that every face has exactly three angles smaller than~$\pi$, some of the vertices have an incident angle greater than $\pi$ and the outer face is bounded by the convex hull of the point set. Different labelings of the angles of plane Laman graphs have been investigated in~\cite{HORSSSSSW,osss-cpt-2007,ST}.

In order to provide plane Laman graphs with a binary labeling, we face the problem that
there are such graphs, e.g. segments or triangles, for which (G$0$) and the edge rule
(G$2$) cannot be simultaneously satisfied, since the two special vertices must be adjacent and the edge
between them would violate (G$2$). Therefore, we have to allow this one
exception and modify the definition of weak labeling:

\begin{definition}\rm
An {\emph{extended weak labeling}} for a plane Laman graph is a
mapping from its angles to the set $\{0,1\}$ which satisfies (G$0'$),
(G1), (G$2'$) and (G3), where

\Item(G$0'$)
{{\bf{Laman special vertices:}} There are two adjacent special
vertices $s_0$ and $s_1$ such that all angles incident to $s_i$ are
labeled $i$.}

\Item(G$2'$) {{\bf{Laman edge rule:}} For each edge except $s_0s_1$, the incident labels
  coincide at one endpoint and differ at the other.}
\end{definition}

\begin{figure}[htb]
    \begin{center}
        \includegraphics{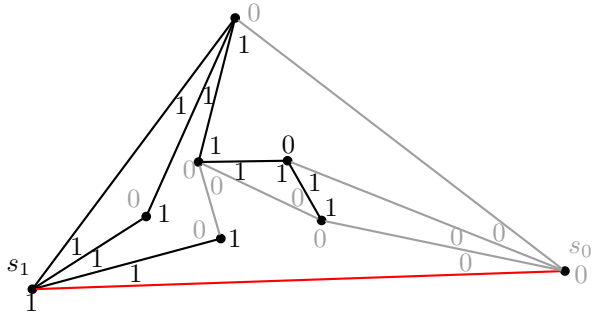}
    \end{center}
    \caption{An extended weak labeling for a pointed pseudo-triangulation.}
    \label{fig:labelling_ppt}
\end{figure}

Figure~\ref{fig:labelling_ppt} shows a extended weak labeling for a plane Laman
graph, embedded as a pointed pseudo-triangulation (see~\cite{Streinu}).

\begin{figure}[htb]
    \begin{center}
        \includegraphics{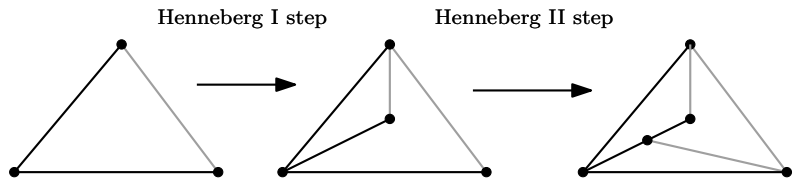}
    \end{center}
    \caption{Constructing a decomposition into two trees via Henneberg steps.}
    \label{fig:henneberg_trees}
\end{figure}

In Theorem~\ref{thm:ppt_labeling} below we show that each plane
Laman graph admits an extended weak labeling. For the proof, we use the
following characterization of Laman graphs: They can be built via
{\emph{Henneberg constructions}}~\cite{HORSSSSSW,osss-cpt-2007,Streinu,TW}, which
start with a triangle and iterate vertex insertions of the following types (see
Figure~\ref{fig:henneberg_trees}):
\begin{itemize}
\item{Add a degree-two vertex (Henneberg I step)}
\item{Place a vertex on an existing edge and connect it to a third vertex (Henneberg II step).}
\end{itemize}

Before proving the existence of extended weak labelings for plane Laman graphs, we need
the following technical result:

\begin{lemma}
\label{lem:topological_Henneberg}
For every plane Laman graph $G$, there exists a Henneberg construction
such that all intermediate graphs are plane and, at each step, the
topology is changed only on edges and faces involved in the
Henneberg step. Furthermore, there exists an edge $e$ of the initial
triangle which is never split in the construction.
\end{lemma}

\begin{proof}
The first part of the statement is Lemma~7 in~\cite{HORSSSSSW}. In order to add the condition in the last sentence, we just have to show that there exists a vertex $v \in G$ of degree $2$ or $3$ which is not an endpoint of~$e$, and then follow the proof there. For proving the claim let us assume, for the sake of a contradiction, that no such vertex $v$ exists. Since $G$ has $2n-3$ edges, the degree sum of all vertices but the endpoints of $e$ is at least $4(n-2)=4n-8$. This implies that the endpoints of $e$ can have degree at most $1$ contradicting $G$ being a Laman graph.
\end{proof}

\begin{theorem}\label{thm:ppt_labeling}
Every plane Laman graph admits an extended weak labeling.
\end{theorem}

\begin{proof}
Starting with a triangle labeled as in
Figure~\ref{fig:henneberg_basis_induction}, Lemma~\ref{lem:topological_Henneberg} allows the extended weak labeling to be maintained at each step of the Henneberg construction:

\begin{figure}[htb]
    \begin{center}
        \includegraphics{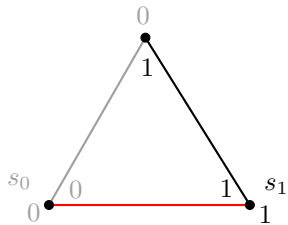}
    \end{center}
    \caption{The initial extended weak labeling in the Henneberg construction.}
    \label{fig:henneberg_basis_induction}
\end{figure}

A Henneberg I step involves only one face of the graph. A new vertex is
placed inside the face and connected to two vertices on the
boundary. Two cases arise depending on whether the two angles at the boundary vertices are
labeled equally or differently: The corresponding
completions of the extended weak labeling are shown in
Figure~\ref{fig:henne1step}. We write $\bar{i}$ for labels $1-i$ in Figures~\ref{fig:henne1step} and~\ref{fig:henne2step}.

\begin{figure}[htb]
    \begin{center}
        \includegraphics{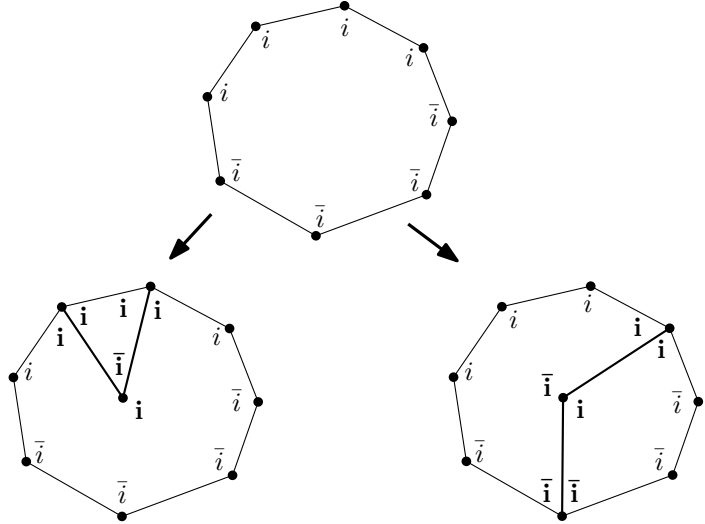}
    \end{center}
    \caption{A Henneberg I step maintains the extended weak labeling.}
    \label{fig:henne1step}
\end{figure}

A Henneberg II step subdivides an edge $e$ and splits one of the two
faces incident to $e$. These two faces are different with respect to
$e$ according to the edge-rule (G2): In one face, the two labels at
$e$ are different, in the other one, the two labels at $e$ are the
same. If we split the face where both angles at $e$ are labeled $i$,
we distinguish two subcases: Either we connect the subdivision vertex
to a vertex with label $i$ or to one with label $1-i$.

\begin{figure}[htb]
    \begin{center}
        \includegraphics[scale=0.7]{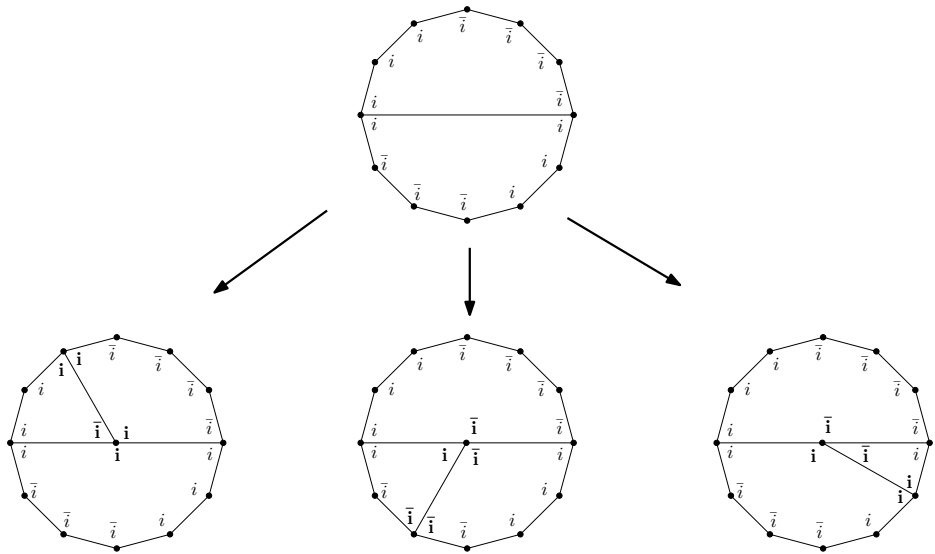}
    \end{center}
    \caption{A Henneberg II step maintains the extended weak labeling.}
    \label{fig:henne2step}
\end{figure}

All three cases and the respective completions of the extended weak labeling
are  shown in Figure~\ref{fig:henne2step}.
\end{proof}

Let us finally note that it is well known that a Laman graph can be decomposed into two trees~\cite{TW}.
These trees can be obtained via the Henneberg
construction, as indicated in Figure~\ref{fig:henneberg_trees}.  The
new vertex is a leaf either in both trees (Henneberg I step) or in
one tree (Henneberg II step).

\begin{figure}[htb]
    \begin{center}
        \includegraphics{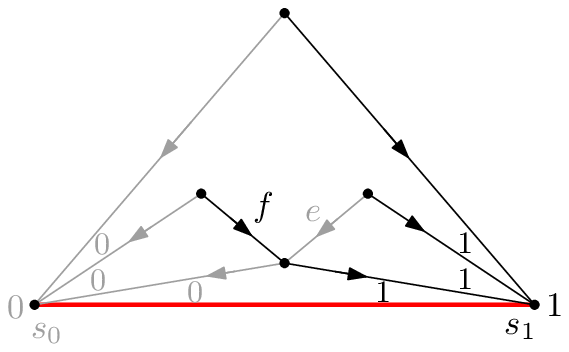}
    \end{center}
    \caption{The extended weak labeling for plane Laman graphs does not induce a decomposition into two trees.}
    \label{fig:no_spanning_trees}
\end{figure}

Unfortunately, although the extended weak labeling is based on the Henneberg
construction too, it does not always give a decomposition of the graph
into two trees; see Figure~\ref{fig:no_spanning_trees} for a simple
example in which the angle formed by~$e$ and~$f$
would have to receive both labels~$0$ and~$1$, contradicting the definition
of extended weak labeling. 

\section*{Acknowledgements}

Apart from the authors' universities, parts of this work were done
during the III Taller de Geometr\'{\i}a Computacional, organized by
the Universidad de Valladolid, and during a visit to the Centre de
Recerca Matem\`{a}tica.

\end{document}